\documentstyle[12pt]{article}
\begin{document}

\title{Holomorphic quantization formula in singular reduction}

\author{Weiping Zhang\thanks{Partially supported by  the NNSF,
SEC of China and the Qiu Shi Foundation.}}

\date{}

\maketitle

\begin{abstract} We show that the holomorphic Morse inequalities proved
by Tian and the author [TZ1, 2] are in effect equalities by refining
the analytic arguments in [TZ1, 2]. 
\end{abstract}

{\bf \S 0. Introduction and the statement of main results}

$\ $

Let $(M,\omega, J)$ be a compact K\"ahler manifold with the K\"ahler form $\omega$ and
the complex structure $J$. 
Let $g^{TM}$ denote the corresponding K\"ahler metric.
We make the assumption that there exists a Hermitian line
bundle $L$ over $M$ admitting a Hermitian connection $\nabla^L$ such that
${\sqrt{-1}\over 2\pi}(\nabla^L)^2=\omega$. 
Then $L$ admits a unique holomorphic 
structure so that $\nabla^L$ is the associated Hermitian holomorphic connection.
We call $L$ the prequantum line bundle over $M$. 

Next, suppose that $(M,\omega, J)$ admits a holomorphic Hamiltonian action 
of a compact connected Lie group $G$ with Lie algebra ${\bf g}$. Let $\mu:M\rightarrow
{\bf g}^*$ be the corresponding moment map. Then a formula due to Kostant [K]
(cf. [TZ1, (1.13)]) induces a natural ${\bf g}$ action on $L$. We make the assumption 
that this ${\bf g}$ action can be lifted to a 
holomorphic $G$ action on $L$. Then this $G$ action
preserves $\nabla^L$. After an integration if necessary, we also assume that this 
$G$ action preserves the Hermitian metric on $L$.
Then for any integer $p\geq 0$, the $p$-th
Dolbeault cohomology $H^{0,p}(M,L)$ is a $G$-representation. We denote its
$G$-invariant part  by $H^{0,p}(M,L)^G$. 

Now let $a\in {\bf g}^*$ be a regular value of $\mu$. Let
${\cal O}_a\subset {\bf g}^*$ be the coadjoint orbit of $a$.
For simplicity,
we assume that $G$ acts on $\mu^{-1}({\cal O}_a)$ freely. Then the quotient 
space $M_{G,a}=\mu^{-1}({\cal O}_a)/G$ is smooth. 
Also, $\omega$ descends canonically
to a symplectic form $\omega_{G,a}$ on $M_{G,a}$ so that one gets the
Marsden-Weinstein reduction $(M_{G,a},\omega_{G,a} )$. Furthermore,
the  complex structure $J$ 
descends canonically to a complex structure $J_{G,a} $ on $TM_{G,a}$ so that
$(M_{G,a}, \omega_{G,a},J_{G,a}) $ is again K\"ahler.
On the other hand, the pair $(L,\nabla^L)$
also descends canonically to a Hermitian holomorphic line bundle $L_{G,a}$ over
$M$. 

One of the purposes of this paper is to present an analytic proof of the 
following result.

$\ $

{\bf Theorem 0.1.} {\it If $\mu^{-1}(0)$ is not empty, then
there exists an open neighborhood ${\bf O}$ of $0\in {\bf g}^*$ such that
for any regular value  $a\in {\bf O}$ of $\mu$ with $\mu^{-1}(a)$
nonempty,  the following identity for Dolbeault 
cohomologies  holds for any integer $p\geq 0$,}
$$\dim H^{0,p}(M,L)^G = \dim H^{0,p}(M_{G,a},L_{G,a}) .\eqno(0.1)$$

In the case where $a=0\in {\bf g}^*$ is a regular value of $\mu$, 
(0.1) was proved by Guillemin-Sternberg [GS] for $p=0$, and recently 
for general $p$ by 
Teleman [T] in the algebraic situation,
 refining in this case the Morse type inequalities 
of Tian and Zhang [TZ1, Theorem 0.4]. 
In this case, (0.1) has also
been extended by Wu [W] to  cases where $M$ is non-compact, when the Lie
group $G$ is Abelian.

While for the case where $0\in {\bf g}^*$ is a singular value of $\mu$, 
Theorem 0.1, of which a more general version has been proved by
Teleman [T, Sect. 5] in the algebraic situation, 
refines the Morse type inequalities proved in [TZ2, Theorem 0.1]. 

In a recent preprint [Br], Braverman proposed an analytic proof of Teleman's
result, in the regular reduction case,
 by extending the methods developed in [TZ1]. In particular, he pointed
out that to get the equality (0.1) for the $a=0$ case, one needs to find 
a quasi-homomorphism between the Dolbeault complexes under considerations
which verifies certain specific properties. Braverman actually
constructed  a quasi-homomorphism in his paper.

We combine the methods and results in [TZ1, 2] with Braverman's idea, that one needs to
construct a suitable quasi-homomorphism,  to prove  Theorem 0.1.
However, the quasi-homomorphism 
we will construct is  different from the one  in [Br],
and ties  closer to the methods in [TZ1].
We first prove an extension of (0.1), when $a=0$ is a regular value of $\mu$,
to the case where the prequantum line bundle $L$ can be replaced by some
more general vector bundles. We then prove Theorem 0.1 by using the 
methods in [TZ2].

The above arguments also lead to the following extension
to the non-Abelian group action case of a result of Wu and Zhang [WZ, Corollary 4.8].

$\ $

{\bf Theorem 0.2.} {\it If $a\in {\bf g}^*$ is a regular value of $\mu$ and
$\mu^{-1}(a)\neq \emptyset$, then
for any integer $p\geq 0$, the following identity holds,}
$$\dim H^{0,p}(M,{\bf C}) = \dim H^{0,p}(M_{G,a},{\bf C}_{G,a}) .\eqno(0.2)$$

This paper is organized as follows. In Section 1, we prove an extended 
holomorphic quantization formula, in the regular reduction case, to allow more  general
coefficients.  In Section 2, we prove Theorems
0.1 and 0.2.

$\ $

{\bf \S 1. An extended holomorphic quantization formula}

$\ $

In this section, we extend the holomorphic quantization formula in [T] 
to allow more general coefficients. Our method is a
 combination of the method in [TZ1, 2] with an idea of Braverman 
[Br] that one needs to construct a quasi-homomorphism between the Dolbeault 
complexes under considerations. 

This section is organized as follows. In a), we state the main result of this
section. In b), we recall the construction of an isometric embedding $J_T$ which
has been essential to the arguments in [TZ1, 2]. In c), we recall some results
in [TZ2] which include in particular Morse type inequalities valid for
coefficients more general than the prequantum line bundle $L$.
In d), we construct the above mentioned quasi-homomorphism. In e), we first prove
an estimate verified by the   quasi-isomorphism constructed in d), which can be 
thought 
of as an analogue of a result of Bismut and Lebeau [BL, (10.4)] in our situation.
We then apply a trick of Braverman [Br] to complete the proof of the main result
stated in a).

$\ $

{\bf a). An extended holomorphic quantization formula}

$\ $

In this section, we assume that $0\in {\bf g}^*$ is a regular value of the moment map
$\mu$ with $\mu^{-1}(0)\neq \emptyset$
 and, for simplicity, that $G$ acts on $\mu^{-1}(0)$ freely. To simplify the notation,
we will denote  $M_{G, 0}$, $\omega_{G, 0}$ etc. 
by $M_{G}$, $\omega_{G}$, etc..

Let $E$ be a $G$-equivariant Hermitian holomorphic vector bundle over $M$ with
the Hermitian holomorphic connection denoted by $\nabla^E$. Then it induces
canonically (cf.  [TZ1]) a Hermitian holomorphic vector bundle $E_{G}$
over $M_{G}$.

Let ${\bf g}$ (and thus ${\bf g}^*$ also) be equipped with an ${\rm Ad}G$-invariant
metric. Let ${\cal H}=|\mu |^2$ be the norm square of the moment map.

Let $h_i$, $1\leq i\leq \dim G$, be an orthonormal base of ${\bf g}^*$. 
Let $V_i$, $1\leq i\leq \dim G$, be the dual base of $h_i$, 
$1\leq i\leq \dim G$.
Then we can write $\mu$ as 
$$\mu=\sum_{i=1}^{\dim G}\mu_ih_i, \eqno (1.1)$$
with each $\mu_i$ a real function on $M$.

For any $V\in {\bf g}$, set\footnote{We use the same notation $V$
to denote the vector field it generates on $M$.}
$$r_V^E=L_V^E-\nabla_V^E,\eqno (1.2)$$
where $L_V^E$ denotes the infinitesimal action of $V$ on $E$.

$\ $

{\bf Assumption 1.1.}  {\it At any critical point $x\in M$
of ${\cal H}$, one has}
$$\sqrt{-1}\sum_{i=1}^{\dim G}\mu_i(x)r_{V_i}^E(x)\geq 0.\eqno (1.3)$$

The main result of this section can be stated as follows,
where we still use the superscript $G$ to denote the $G$-invariant part.

$\ $

{\bf Theorem 1.2.}  {\it If $E$ verifies Assumption 1.1, then 
the following identity for Dolbeault
cohomologies holds for any integer $p\geq 0$,}
$$\dim H^{0,p}(M,E)^G= \dim H^{0,p}(M_{G},E_{G}) .\eqno(1.4)$$

{\bf Remark 1.3.} Clearly, the prequantum line bundle $L$ verifies the Assumption
1.1. Thus (1.4) holds for  $E=L$ and $E={\bf C}$. That is,
it contains  the Teleman formula [T] 
and Theorem 0.2 in the case where $a=0$ is a regular value of $\mu$. 

$\ $

{\bf b). An embedding from $\Omega^{0,*}(M_{G},E_{G}) $ into
$\Omega_G^{0,*}(M,E) $}

$\ $

Let $\Omega^{0,*}(M,E) $ (resp. $\Omega^{0,*}(M_{G},E_{G}) $)
denote the set of smooth sections of $\wedge^{0,*}(T^*M)\otimes E $ 
(resp. $\wedge^{0,*}(T^*M_{G})\otimes E_{G} $). We denote by
$\Omega_G^{0,*}(M,E) $ the $G$-invariant part of 
$\Omega^{0,*}(M,E) $. 

In this subsection, following [BL] and
[TZ1], we construct explicitly, for any $T>0$, an isometric embedding 
$J_T:\Omega^{0,*}(M_{G},E_{G}) \rightarrow \Omega_G^{0,*}(M,E) $.

As in [TZ1, Sect. 3b)], let $U$ be a sufficiently small $G$-invariant open
neighborhood of $\mu^{-1}(0)$ such that $G$  acts on $U$ freely. 
Then $U/G$ is  smooth  and carries a canonically induced 
metric $g^{T(U/G)}$ (cf. [TZ1, (3.8)]). 

Let $N_G$ be the normal bundle to 
$M_{G} $ in $U/G$. We identify $N_G$ as the orthogonal completement of 
$TM_{G} $ in $ T(U/G)$, i.e., 
$T(U/G)|_{M_{G} } =N_G\oplus TM_{G} $ and
$g^{T(U/G)}|_{M_{G}}=g^{N_G}\oplus g^{TM_{G} },$ 
where $g^{N_G}$ is the induced metric on $N_G$.

If $y\in M_{G} $, $Z\in N_{G,y}$, let 
$t\in {\bf R} \mapsto x_t=\exp_y^{U/G}(tZ) \in U/G$
be the geodesic in $U/G$ with $x_0=y$, ${dx_t\over dt}|_{t=0}=Z$. For $ \varepsilon>0$,
set $B_{\varepsilon } =\{ Z\in N_G;|Z|<\varepsilon \}$. Since $M_{G} $ is
compact, there exists $\varepsilon_0 >0$ such that for 
$0<\varepsilon <\varepsilon_0 $, the map 
$(y,Z)\in N_G\mapsto \exp_y^{U/G}(Z) \in U/G$ is 
a diffeomorphism from $B_{\varepsilon } $ onto a tubular neighborhood 
$U_{\varepsilon } $ of $M_{G} $ in $U/G$. From now on, we identify 
$B_{\varepsilon } $ with $U_{\varepsilon } $ and use the notation $(y,Z)$
instead of $\exp_y^{U/G}(Z) $. In particular, we identify $y\in M_{G} $ with
$(y,0)\in N_G$.

Let $dv_{N_G}$ be the volume form of the fibers in $N_G$. Then 
$dv_{M_{G} }(y) dv_{N_G}(Z) $ is a natural volume form on the total space of $N_G$.
Let $k(y,Z)$ be the smooth positive function on $B_{\varepsilon_0 } $ 
defined by  $dv_{U/G}(y,Z)=k(y,Z)dv_{M_{G} }(y) dv_{N_G}(Z) $. 
The function $k$
has a positive lower bound on $B_{\varepsilon_0 /2} $. Also, $k(y)=1$. 

Now consider the fibration
$G \rightarrow U\stackrel{\pi}{\rightarrow} U/G$. 
Let $h$ be the positive function
on $U/G$ defined by
${h}(x)=\sqrt{{\rm vol} (\pi^{-1}(x))}$ for any $x\in U/G$.

If $F$ is a $G$-invariant Hermitian vector bundle over $U$, then it induces
canonically a Hermitian vector bundle $F_{U/G} $ over $U/G$ such that 
$\pi^*F_{U/G}=F$. We  denote as in [TZ1, Sect. 3c)] by $\pi_G^F :\Gamma_G(F)
\rightarrow \Gamma (F_{U/G} )$ the canonical isomorphism which maps a
$G$-invariant section of $F$ to the corresponding section of $F_{U/G} $. 
We will usually omit the superscript $F$ from $\pi_G^F $ if there will be
no confusion. 
If $\nabla^F$ is a $G$-invariant 
Hermitian connection on $F$, then it induces canonically a Hermitian connection
$\nabla^{F_{U/G}}$ on $F_{U/G}$ (cf. [TZ1, (3.9)]). 

Now as in [BL, Sect. 8g)], for $x=(y,Z)\in U_{\varepsilon_0 } $, we identify
$E_{U/G,x}$ (resp. $(\wedge^{0,*}(T^*M))_{U/G,x}$) 
with $E_{U/G,y}$ (resp. $(\wedge^{0,*}(T^*M))_{U/G,y}$) 
by parallel transport with respect to $\nabla^{E_{U/G}}$ 
(resp. $\nabla^{(\wedge^{0,*}(T^*M))_{U/G}}$, with
 $\nabla^{\wedge^{0,*}(T^*M) }$ the Hermitian holomorphic connection on
$\wedge^{0,*}(T^*M) $) along the geodesic 
$t\mapsto (y,tZ)$. The induced identification of 
$(\wedge^{0,*}(T^*M)\otimes E)_{U/G,x}$ 
with $(\wedge^{0,*}(T^*M)\otimes E)_{U/G,y}$ preserves the  
the ${\bf Z}$-grading
of $(\wedge^{0,*}(T^*M)\otimes E)_{U/G}$ and is  $G$-equivariant.

Since $0\in {\bf g}^*$ is a regular value of the moment map 
$\mu:M\rightarrow {\bf g}^*$, one verifies easily that $\mu^{-1}(0)$
is a nondegenerate submanifold, in the sense of Bott, of 
${\cal H}=|\mu |^2$.
Let ${d}\mu:TM|_{\mu^{-1}(0)} \rightarrow {\bf g}^*$ denote 
the restriction of the differential of $\mu$ on $\mu^{-1}(0)$.
Clearly, for any $Z\in N_G$, $\pi^*Z\in TM|_{\mu^{-1}(0)} $ and $|{d}\mu(\pi^*Z)|$
depends only on $Z$.

Take $\varepsilon \in (0,{\varepsilon _0/2}]$. Let 
$\rho :{\bf R}\rightarrow [0,1]$ be a smooth function such that
$\rho (a)=1$ if $a\leq 1/2$, while $\rho (a)=0$ if $a\geq 1$. 
For $Z\in N_G$, set $\rho_{\varepsilon }(Z)=\rho (|Z|/\varepsilon )$.

For  $T>0$, $y\in M_{G}$, set 
$$\alpha_T(y)=\int_{N_{G,y}}\exp \left(-T|{d}\mu(\pi^*Z)|^2\right) 
\rho_{\varepsilon }^2(Z) dv_{N_G}(Z).\eqno (1.5)$$

$\ $

{\bf Definition 1.4.} {\it For any $T>0$, let 
$J_T:\Omega^{0,*}(M_{G},E_{G}) \rightarrow \Omega^{0,*}_G(M,E) $ be defined by}
$$J_T:\alpha \mapsto \pi^* \left(
k^{-1/2}\alpha_T^{-1/2}h^{-1} \rho_{\varepsilon }(Z)\exp \left(-
{T|{d}\mu(\pi^*Z)|^2\over 2} \right)
\alpha \right).\eqno(1.6)$$

One verifies easily that $J_T$ is well-defined and that it is an isometry from 
$\Omega^{0,*}(M_{G},E_{G}) $ onto its image. In fact, although not explicitly
written out, the map $J_T$ has played an essential role in [TZ1, 2],
in particular, in getting the results in the following subsection.

$\ $

{\bf c). Deformations of  Dolbeault complexes and Morse type inequalities}

$\ $

In the rest of this section, we assume that $E$ verifies Assumption 1.1.

Let $\overline{\partial}^E $ (resp. $\overline{\partial}^{E_{G}} $)
 be the Dolbeault  operator acting on 
$\Omega^{0,*}(M,E) $ (resp.
$\Omega^{0,*}(M_{G},E_{G}) $). 

Following Tian and Zhang [TZ1, (1.21)], for any $T\in {\bf R}$, set
$$\overline{\partial}^E_T =e^{-T|\mu |^2/2} \overline{\partial}^E e^{T|\mu |^2/2}:
\Omega^{0,*}(M,E) \rightarrow \Omega^{0,*}(M,E) ,$$ 
$$ D_T^E = \sqrt{2}\left( \overline{\partial}^E_T +
\left(\overline{\partial}^E_T\right)^* \right) . \eqno (1.7)$$

Set $\widetilde{h}=h|_{M_{G}}$. 
Following [TZ1, (3.54)], set
$$ \overline{\partial}^{E_{G}}_Q =\widetilde{h} \overline{\partial}^{E_{G}} \widetilde{h}^{-1}: 
\Omega^{0,*}(M_{G},E_{G}) \rightarrow \Omega^{0,*}(M_{G},E_{G})  ,$$
$$D^{E_{G}}_Q= \sqrt{2}\left( \overline{\partial}^{E_{G}}_Q +
\left( \overline{\partial}^{E_{G}}_Q\right)^*\right)  . \eqno (1.8)$$

By proceeding as in [TZ1, Sects. 2-4], 
one gets easily the  following  refinement of [TZ1, Theorems 3.13 (ii)
and 4.2].

$\ $

{\bf Proposition 1.5.} {\it There exist $c_0>0$, $T_0>0$ such that
there are no nonzero eigenvalues of $D_Q^{E_{G},2}$ in $[0,c_0]$, and that
for any $T\geq T_0$, the number of eigenvalues of 
$D_T^{E,2}|_{ \Omega_G^{0,*}(M,E) }$ in $[0,c_0]$ equals to 
$\dim (\ker D_Q^{E_{G}})$.}

$\ $

{}From Proposition 1.5 and the ${\bf Z}$-grading nature of the problem, one gets for
any integer $p\geq 0$ the following Morse type inequalities proved in 
[TZ2, Theorem 2.1],
$$ \dim H^{0,p}(M,E)^G \leq \dim H^{0,p}(M_{G},E_{G}).\eqno (1.9)$$

$\ $

{\bf d). A quasi-homomorphism $r:(\Omega_G^{0,*}(M,E),\overline{\partial}^E_T ) \rightarrow 
(\Omega^{0,*}(M_{G},E_{G}), \overline{\partial}^{E_{G}}_Q )$ }

$\ $

In this subsection, we construct a  quasi-homomorphism from
$ (\Omega_G^{0,*}(M,E),\overline{\partial}^E_T ) $ to
$(\Omega^{0,*}(M_{G},E_{G}), \overline{\partial}^{E_{G}}_Q ) $.
 Our quasi-homomorphism is different from 
the one constructed by Braverman in [Br].

Let $i:\mu^{-1}(0)\hookrightarrow M$ denote the canonical isometric embedding.

One  verifies easily that the induced Hermitian holomorphic bundle $E_G$ from
$E$ is given by $E_G=(i^*E)_{U/G} $ 
with its Hermitian holomorphic connection
given by $\nabla^{(i^*E)_{U/G}} $. 

Let $N=\pi^*N_G$ be the normal bundle to $\mu^{-1}(0) $ in $M$. Then one verifies
that $JN$ is the vertical tangent vector bundle of the fibration
$G \rightarrow \mu^{-1}(0) \stackrel{\pi}{\rightarrow} M_G $ (cf. [TZ1, Sect. 3]). 
Thus one has the canonical orthogonal splittings
$$TM|_{\mu^{-1}(0) } =N\oplus JN \oplus \pi^*(TM_G) ,$$
$$g^{TM}|_{\mu^{-1}(0) } =g^N\oplus g^{JN} \oplus \pi^*g^{TM_G},\eqno (1.10)$$
with $J$ preserves $N_J=N\oplus JN $ and $W_J=\pi^*(TM_G) $. Thus one has as in
[TZ1, (3.40)] the canonical identifications of Hermitian vector bundles
$$\wedge^{0,*} (T^*M)|_{\mu^{-1}(0) } = \wedge^{0,*} (N_{J}^* )\widehat{\otimes}
\wedge^{0,*} (W_J^* ) .\eqno (1.11)$$

Let $q$ be the canonical orthogonal projection
$$q: \wedge^{0,*} (N_{J}^* ) \widehat{\otimes}
\wedge^{0,*} (W_J^* ) \otimes i^*E \rightarrow \wedge^{0,*} (W_J^* ) \otimes i^*E ,
\eqno (1.12)$$
which acts as identity on $\wedge^{0,0} (N_{J}^* ){\otimes}
\wedge^{0,*} (W_J^* ) \otimes i^*E \simeq \wedge^{0,*} (W_J^* ) \otimes i^*E  $
and maps each $\wedge^{0,i} (N_{J}^* ) {\otimes}
\wedge^{0,*} (W_J^* ) \otimes i^*E $, $i\geq 1$, to zero.

$\ $

{\bf Proposition 1.6.} {\it The following identity  holds,}
$$  \pi_G qi^* \overline{\partial}^E =\overline{\partial}^{E_{G}} \pi_G qi^* :
\Omega_G^{0,*}(M,E) \rightarrow 
\Omega^{0,*}(M_{G},E_{G}) .\eqno (1.13)$$

{\it Proof.} For any $e\in TM$, write its complexification as $e=e^{1,0}+e^{1,0}$
with $e^{1,0} \in T^{(1,0)}M$, $e^{0,1}\in T^{(0,1)}M$. Let $c(e)$ be the 
Clifford action on $\wedge^{0,*}(T^*M)\otimes E$ defined by 
$$c(e)=\sqrt{2} \left(\overline{e^{1,0}}^*\wedge -i_{e^{0,1}}\right) ,\eqno(1.14)$$
where $\overline{e^{1,0}}^* \in T^{(0,1)*}M$ is the metric dual of 
$e^{1,0}$. Then one verifies easily that 
$$c(Je)=-\sqrt{-2}\left(\overline{e^{1,0}}^*\wedge +i_{e^{0,1}}\right) .\eqno(1.15)$$

Let $f_1,\cdots,f_{\dim M}$ be an 
orthonormal base of $TM$. 
Let $\nabla^{\wedge^{0,*}(T^*M)\otimes E } $ be the Hermitian holomorphic 
connection on 
$\wedge^{0,*}(T^*M)\otimes E $. Then one verifies easily that 
$$\overline{\partial}^E =\sum_{i=1}^{\dim M} \overline{f_i^{1,0}}^*\wedge
\nabla^{\wedge^{0,*}(T^*M)\otimes E }_{f_i} ,\eqno (1.16)$$
by which
one finds that if $e_1,\cdots,e_{\dim M_G}$ is 
an orthonormal base of $W_J $, then
$$qi^* \overline{\partial}^E =\sum_{i=1}^{\dim M_G} \overline{e_i^{1,0}}^*\wedge
\left(i^* \nabla^{\wedge^{0,*}(T^*M)\otimes E }({e_i})\right) i^* .\eqno (1.17)$$

Now, let $P$ (resp. $P^\perp$) be the orthogonal projection from 
$TM|_{\mu^{-1}(0)}$ to $N_J$ (resp. $W_J$) with respect to the orthogonal 
splitting (1.10). Set
$$\nabla^{N_J}=Pi^*\nabla^{TM}P ,\ \ \nabla^{W_J}=P^\perp i^*\nabla^{TM}P^\perp,
\eqno (1.18)$$
where $\nabla^{TM}$ is the Levi-Civita connection of $g^{TM}$, and
$$A=i^*\nabla^{TM}-\nabla^{N_J}-\nabla^{W_J}.\eqno (1.19)$$
Then $A$ and $J$ commute with each other.

As in [TZ1, pp. 252], $\nabla^{N_J}$, $\nabla^{W_J}$ induce canonically
the Hermitian connections $\nabla^{\wedge^{0,*}(N_J^*) }$,  
$\nabla^{\wedge^{0,*}(W_J^*) }$ on $\wedge^{0,*}(N_J^*) $,
$\wedge^{0,*}(W_J^*)$. Let 
$^0\nabla^{(\wedge^{0,*}(T^*M)\otimes E)|_{\mu^{-1}(0)}}$ 
denote the tensor product connection of  $\nabla^{\wedge^{0,*}(N_J^*) }$,  
$\nabla^{\wedge^{0,*}(W_J^*) }$ and $i^*\nabla^E$. Then by [TZ1, (3.46)], one 
has that if $e_{\dim M_G+1},\cdots,e_{\dim \mu^{-1}(0)}$ is an orthonormal
base of $JN$, then for any
$1\leq i\leq \dim M_G$, 
$$i^*\nabla^{\wedge^{0,*}(T^*M)\otimes E } ({e_i}) 
=\, ^0\nabla^{(\wedge^{0,*}(T^*M)\otimes E)|_{\mu^{-1}(0)}}_{e_i} $$
$$+{1\over 2}\sum_{s=1}^{\dim M_G } \sum_{t=\dim M_G +1}^{\dim \mu^{-1}(0)} 
\left(\langle A(e_i)e_s,e_t\rangle c(e_s)c(e_t)
+  \langle A(e_i)e_s,Je_t\rangle c(e_s)c(Je_t)\right).\eqno (1.20)$$

Now one verifies directly that
$$\pi_G q \sum_{i=1}^{\dim M_G} \overline{e_i^{1,0}}^* \wedge \,\left(
^0\nabla^{(\wedge^{0,*}(T^*M)\otimes E)|_{\mu^{-1}(0)}}_{e_i}\right) i^* =
\overline{\partial}^{E_{G}} \pi_G qi^* .\eqno (1.21)$$

On the other hand, by (1.14), (1.15) and the fact that $A$ and $J$ commute
with each other,  one deduces that
$${1\over 2} \sum_{s=1}^{\dim M_G } \sum_{t=\dim M_G +1}^{\dim \mu^{-1}(0)} 
\left(\langle A(e_i)e_s,e_t\rangle c(e_s)c(e_t)
+  \langle A(e_i)e_s,Je_t\rangle c(e_s)c(Je_t)\right)$$
$$={1\over 2}\sum_{s=1}^{\dim M_G } \sum_{t=\dim M_G +1}^{\dim \mu^{-1}(0)} 
\left(\langle A(e_i)e_s,e_t\rangle c(e_s)c(e_t)
+  \langle A(e_i)e_s,e_t\rangle c(Je_s)c(Je_t)\right)$$
$$=- 2\sum_{s=1}^{\dim M_G } \sum_{t=\dim M_G +1}^{\dim \mu^{-1}(0)} 
\langle A(e_i)e_s,e_t\rangle  \left( \overline{e_s^{1,0}}^* \wedge 
i_{e_t^{0,1}} + i_{e_s^{0,1}} \overline{e_t^{1,0}}^* \wedge \right).\eqno (1.22)$$

{\bf Lemma 1.7.} {\it For any $\dim M_G+1\leq t \leq \dim \mu^{-1}(0)$, 
one has,}
$$\sum_{i=1}^{\dim M_G } \sum_{s=1}^{\dim M_G } \langle A(e_i)e_s,e_t \rangle  
\overline{e_i^{1,0}}^* \wedge \overline{e_s^{1,0}}^* =0.\eqno (1.23)$$

{\it Proof.}  By the fibration structure of 
$G \rightarrow \mu^{-1}(0) \stackrel{\pi}{\rightarrow} M_G $ , one knows that
$[e_i, e_t]$'s and $[Je_i, e_t]$'s, with
$1\leq i\leq \dim M_G$, $\dim M_G+1\leq t\leq \dim \mu^{-1}(0)$,
 are lying in $\Gamma(JN)$. Thus one has, since $A$ and $J$
commute with each other, 
$$ \langle A(e_i)e_s,e_t \rangle =\langle \nabla_{e_i}e_s,e_t \rangle 
= -\langle e_s,\nabla_{e_i} e_t \rangle 
= -\langle e_s,\nabla_{e_t} e_i \rangle$$
$$ = -\langle Je_s,\nabla_{e_t} Je_i \rangle 
=\langle A(Je_i)Je_s,e_t \rangle .\eqno (1.24)$$

By (1.24), (1.14) and (1.15), one deduces that
$$\sum_{i=1}^{\dim M_G } \sum_{s=1}^{\dim M_G } \langle A(e_i)e_s,e_t \rangle  
\overline{e_i^{1,0}}^* \wedge \overline{e_s^{1,0}}^* $$
$$=\sum_{i=1}^{\dim M_G } \sum_{s=1}^{\dim M_G } \langle A(Je_i)Je_s,e_t \rangle  
\overline{(Je_i)^{1,0}}^* \wedge \overline{(Je_s)^{1,0}}^* $$
$$=-\sum_{i=1}^{\dim M_G } \sum_{s=1}^{\dim M_G } \langle A(e_i)e_s,e_t \rangle  
\overline{e_i^{1,0}}^* \wedge \overline{e_s^{1,0}}^*, \eqno (1.25)$$
from which (1.23) follows.  $\Box$

$\ $

{}From (1.23) and (1.22), one finds that
$$q  \sum_{i=1}^{\dim M_G } \sum_{s=1}^{\dim M_G } 
\sum_{t=\dim M_G +1}^{\dim \mu^{-1}(0)} \left(
\langle A(e_i)e_s,e_t\rangle \overline{e_i^{1,0}}^*
\wedge c(e_s)c(e_t)\right. $$
$$+\left.
\langle A(e_i)e_s,Je_t\rangle \overline{e_i^{1,0}}^*\wedge c(e_s)c(Je_t) 
\right) i^*=0. \eqno(1.26)$$

{}From (1.17), (1.20), (1.21) and (1.26), one gets (1.13). $\Box$

$\ $

{\bf Definition 1.8.} {\it Let $r$ be the map defined by}
$$r=\widetilde{h} \pi_G qi^* : \Omega_G^{0,*}(M,E) \rightarrow 
\Omega^{0,*}(M_{G},E_{G}) .\eqno(1.27)$$ 

{\bf Theorem 1.9.} {\it For any $T\in {\bf R}$, the map
$$r: \left(\Omega_G^{0,*}(M,E),\overline{\partial}^E_T \right) \rightarrow 
\left(\Omega^{0,*}(M_{G},E_{G}), \overline{\partial}^{E_{G}}_Q \right) \eqno(1.28)$$
is a quasi-homomorphism.}

{\it Proof.} Let $X^{\cal H} $ be the Hamiltonian vector field associated with
${\cal H}=|\mu|^2$. By (1.7) one deduces that 
$$\overline{\partial}^E_T =\overline{\partial}^E + {\sqrt{-1}T\over 2}
\overline{\left( X^{\cal H} \right)^{1,0}}^* \wedge .\eqno (1.29)$$

On the other hand, one verifies easily that
$$\left. X^{\cal H}\right|_{\mu^{-1}(0)}=0.\eqno(1.30)$$

{}From (1.8), (1.13), (1.27), (1.29) and (1.30), one gets
$$r\overline{\partial}^E_T =\overline{\partial}^{E_{G}}_Q r:
\Omega_G^{0,*}(M,E) \rightarrow 
\Omega^{0,*}(M_{G},E_{G}) ,\eqno (1.31)$$
from which Theorem 1.9 follows.  $\Box$

$\ $

{\bf e). Proof of Theorem 1.2}

$\ $

Recall that $c_0>0$ and $T_0>0$ have been fixed in Proposition 1.5.
For any $T\geq T_0$, let $F^{[0,c_0]}_{G,T} $ be the finite dimensional subspace of
$\Omega_G^{0,*}(M,E) $ generated by the eigenspaces of those eigenvalues of
$D_T^{E,2}|_{ \Omega_G^{0,*}(M,E) }$ which lie in $[0,c_0]$.
Let $P_T$ denote the orthogonal projection from 
$\Omega_G^{0,*}(M,E) $ onto  $F^{[0,c_0]}_{G,T}$.

Let $P_Q$ denote the orthogonal projection from $\Omega^{0,*}(M_{G},E_{G}) $ to
$\ker (D^{E_G}_Q) $.

We first prove the following analogue of [BL, (10.4)].

$\ $

{\bf Theorem 1.10.} {\it There exist $\varepsilon >0$,
$C>0$, $T_1>0$ such that for any
$T\geq T_1$, any $\sigma \in \ker (D^{E_G}_Q) $,}
$$\left\| P_Q rP_T\sqrt{\alpha_T} J_T\sigma -\sigma \right\|_0 \leq 
{C\over \sqrt{T}} \| \sigma \|_0 .\eqno(1.32)$$

{\it Proof.}  We choose $T_0$ in Proposition 1.5 so that $c_0$ is not
an eigenvalue of $D_T^{E,2}|_{ \Omega_G^{0,*}(M,E) } $, $T\geq T_0$. 
Let $\delta$ be the circle of center $0$ with radius $\sqrt{c_0}$
in ${\bf C}$. Then one has the following 
analogue of [BL, (10.6)],
$$P_TJ_T\sigma ={1\over 2\pi \sqrt{-1}} \int_\delta \left( \lambda -
D_T^{E}|_{ \Omega_G^{0,*}(M,E) } \right)^{-1}J_T\sigma d\lambda .\eqno (1.33)$$

One can also show that, by proceeding as in [BL, Sect. 9] and [TZ1, 2],
  for $T$ large
enough and $\lambda \in \delta$, 
$\| (\lambda -D_T^{E}|_{ \Omega_G^{0,*}(M,E) } )^{-1}\|_\infty$ is uniformly
bounded. 

The analogues of the above two facts in [BL, Proof of (10.4)] are all the 
ones in [BL, Proof of (10.4)] which were proved by using [BL, Theorem 9.25]. 
The point now is that here we need not an analogue of [BL, Theorem 9.25]  to
have these two properties. Proposition 1.5 is enough for our purpose.

On the other hand, the analogue of [BL, (10.29)], which now reads
$$r {1\over 2\pi \sqrt{-1}} \int_\delta 
\pi^* \left( h^{-1} \exp \left(- {|d\mu (\pi^*Z)|^2\over 2}\right) 
{\sigma \over \lambda} \right)d\lambda =\sigma ,\eqno(1.34)$$
clearly holds.
One can then proceed as in [BL, Proof of (10.4)] to complete the proof of
Theorem 1.10. $\Box$

$\ $

{\bf Proof of Theorem 1.2.} We proceed as in [Br, Sect. 3].

First of all, from (1.32) one knows that
$$P_T\sqrt{\alpha_T}J_T: \ker\left(D^{E_G}_Q \right) \rightarrow 
F^{[0,c_0]}_{G,T} \eqno(1.35)$$
is  surjective when $T$ is very large. Combining with Proposition 1.5, we see
that it is in fact an isomorphism.

Take $\alpha \in F^{[0,c_0]}_{G,T} $. Then 
$ \overline{\partial}^E_T \alpha \in F^{[0,c_0]}_{G,T} $. By the above discussion,
there exists $\beta \in \ker (D^{E_G}_Q )$ such that
$$ \overline{\partial}^E_T \alpha =P_T\sqrt{\alpha_T}J_T\beta .\eqno(1.36)$$

{}From (1.36) and (1.31), one finds
$$P_QrP_T\sqrt{\alpha_T}J_T\beta  =P_Qr\overline{\partial}^E_T \alpha 
=P_Q\overline{\partial}^{E_G}_Qr \alpha =0.\eqno(1.37)$$

{}From (1.37) and (1.32), one finds
$$\| \beta \|_0 \leq {C\over \sqrt{T}} \| \beta \|_0 ,\eqno(1.38)$$
from which one sees that $\beta=0$ as $T$ is large enough.

Thus, when $T$ is large enough, one has that
$$\left. \overline{\partial}^E_T \right| _{F^{[0,c_0]}_{G,T} }=0.\eqno (1.39)$$

{}From (1.39) and Proposition 1.5, one sees that when $T$ is large enough,
$$\dim \ker \left(D_T^{E}|_{ \Omega_G^{0,*}(M,E) } \right) =\dim
F^{[0,c_0]}_{G,T} =\dim \ker\left(D^{E_G}_Q\right).\eqno(1.40)$$

By (1.40) and (1.9), one completes the proof of Theorem 1.2. $\Box$

$$\ $$

{\bf \S 2. Proof of Theorems 0.1 and 0.2}

$\ $

This section is organized as follows. In a),  we apply Theorem 1.2 
and a trick in [TZ2]
to prove Theorem 0.1. In b), we prove Theorem 0.2.

$\ $

{\bf a). Proof of Theorem 0.1}

$\ $

We proceed as in [TZ2]. Take $a\in {\bf g}^*$.
Since the coadjoint orbit ${\cal O}_a$ admits a 
canonical K\"ahler 
form $\omega_a$ and the holomorphical Ad$ G$ action on 
${\cal O}_a$ is Hamiltonian with the moment map given by the 
canonical embedding $i_a:{\cal O}_a\hookrightarrow {\bf g}^*$ (cf. 
[MS, Chap. 5]),  the induced action
of $G$ on the K\"ahler product 
$(M\times {\cal O}_a,\omega\times (-\omega_a))$ is also 
Hamiltonian with the moment map
${\mu}_a:M\times {\cal O}_a\rightarrow {\bf g}^*$ given by
$${\mu}_a(x,b)=\mu(x)-b.\eqno (2.1)$$

Set ${\cal H}_a=|{\mu}_a|^2$.
The following result has been proved in [TZ2].  

$\ $

{\bf Lemma 2.1.} (Tian-Zhang [TZ2, Prop. 1.2])
{\it There exists an open neighborhood ${\bf O}$ of
$0\in {\bf g}^*$ such that if $a\in {\bf O}$ and $(x,b)\in M\times {\cal O}_a$
is a critical point of ${\cal H}_a$, then the following 
inequality for the inner product on ${\bf g}^*$ holds,}
$$\langle \mu(x)-b,\mu(x)\rangle \geq 0.\eqno (2.2)$$

Now let $a\in {\bf O}$ be a regular value of $\mu$ and,
for simplicity, that $G$ acts on $\mu^{-1}(a)$ freely.
Then $0\in {\bf g}^*$ is a regular value of ${\mu}_a$.
Furthermore,
one has the standard identification of the symplectic 
quotients 
$${\mu}^{-1}_a(0)/G \equiv \mu^{-1}({\cal O}_a)/G
=M_{G,a}\eqno (2.3)$$
carrying the canonically induced K\"ahler form $\omega_{G,a}$.

Recall that $L$ is the holomorphic pre-quantum line bundle over $M$.
Let $\pi$ denote the projection from
$M\times {\cal O}_a$ to its first factor $M$. Let
${\cal L}=\pi^*L$ be the pull-back holomorphic Hermitian
line bundle over  $M\times {\cal O}_a$, with the Hermitian holomorphic 
connection denoted by $\nabla^{\cal L}$.
Then the $G$ action on $L$ lifts canonically
to a holomorphic 
 $G$-action on ${\cal L}$. In particular, for any $V\in {\bf g}$,
its  infinitesimal action on ${\cal L}$ is, via
the Kostant formula [K] (cf. [TZ1, (1.13)]) for the ${\bf g}$-action
on $L$,    given by
$$L^{\cal L}_V=\nabla^{\cal L}_V-2\pi \sqrt{-1}\langle \mu \pi,
V\rangle ,\eqno (2.4)$$
from which one has, in using the notation  in (1.2), that
$$r^{\cal L}_{V_i}=-2\pi \sqrt{-1} \mu_i(x)\eqno (2.5)$$
at any point $(x,b)\in M\times {\cal O}_a$.
By (2.1), (2.5) and Lemma 2.1, one verifies that at any critical point
$(x,b)\in M\times {\cal O}_a$
of ${\cal H}_a$, 
$$\sum_{i=1}^{\dim G} \sqrt{-1}{\mu}_{a,i}(x,b)r^{\cal L}_{V_i}
(x,b)= 2\pi \langle
\mu(x)-b,\mu(x)\rangle \geq 0.\eqno (2.6)$$
 
On the other hand, one verifies directly that the induced
line bundle ${\cal L}_G$ over ${\mu}_a^{-1}(0)/G$ is exactly the line
bundle $L_{G,a}$   over $M_{G,a}=\mu^{-1}({\cal O}_a)/G$.

One can then apply  Theorem 1.2 to $M\times {\cal O}_a$,
${\mu}_a$ and ${\cal L}$ to get
$$\dim H^{0,p}(M\times {\cal O}_a,{\cal L})^G 
=\dim H^{0,p}(M_{G,a},L_{G,a}) \eqno(2.7)$$
for any integer $ p\geq 0$. 

Furthermore, by the definition of ${\cal L}$, one verifies directly 
that
$$\dim H^{0,p}(M\times {\cal O}_a,{\cal L})^G=\sum_{i+j=p}
\dim H^{0,i}(M,L)^G \cdot \dim H^{0,j} ({\cal O}_a,{\bf C})^G$$
$$=\dim H^{0,p}(M,L)^G, \eqno (2.8)$$
with the last equality follows from the easy facts that
$$\dim H^{0,0}({\cal O}_a, {\bf C}) =
\dim H^{0,0}({\cal O}_a, {\bf C}) ^G=1 ,$$
$$\dim H^{0,j}({\cal O}_a, {\bf C})=0 \ \ {\rm for}\ \ j\geq 1.\eqno (2.9)$$

(0.1) follows from (2.7) and (2.8).  $\Box$

$\ $

{\bf b). Proof of theorem 0.2}

$\ $

The main observation is that, in using the notation in (1.2), one has 
$$ r_V^{\bf C} \equiv 0\eqno (2.10)$$
for any $V\in {\bf g}$.
Thus ${\bf C}$ verifies the conditions of Theorem 1.2. Furthermore,
(2.10) plays  roles for ${\bf C}$ similar to what Lemma 2.1 plays for ${\cal L}$
in a).  And 
one sees that to proceed the arguments in a) for ${\bf C}$, one no longer needs to
assume that $a$ is close to $0$ in ${\bf g}^*$. That is, one can show that for
any regular value $a\in {\bf g}^*$ of $\mu$, 
$$\dim H^{0,p}(M,{\bf C})^G =\dim H^{0,p}(M_{G,a},{\bf C}_{G,a}) .\eqno(2.11)$$

On the other hand, since $G$ is connected, one verifies easily that
$$\dim H^{0,p}(M,{\bf C}) =\dim H^{0,p}(M,{\bf C})^G .\eqno (2.12)$$

By (2.11) and (2.12), the proof of Theorem 0.2 is completed. $\Box$

$\ $

{\bf Remark 2.2.} Using the arguments in this paper, one can also show that the 
`weighted' Morse type inequalities proved in [TZ3] are in effect equalities.

$\ $

{\small {\bf Acknowledgements.}  
We are indebted to Youliang Tian in an obvious way for the
many ideas he shared with us through the joint works [TZ1, 2]. }

$\ $

{\bf References}

$\ $

{\small 
[BL] J.-M. Bismut and G. Lebeau, Complex immersions and Quillen metrics. 
{\it Publ. Math. IHES.} 74 (1991), 1-297.

[Br] M. Braverman, Cohomology of the Mumford quotient. {\it Preprint}, \\
math.SG/9809146.

[GS] V. Guillemin and S. Sternberg, Geometric quantization and
multiplicities of group representations. 
{\it Invent. Math.} 67 (1982), 515-538.

[K] B. Kostant, Quantization and unitary representations. {in} {\it Modern
Analysis and Applications.}  Lecture Notes in Math. vol. 170,
Springer-Verlag, (1970), pp. 87-207.

[MS] D. Mcduff and D. Salamon, {\it Introduction to Symplectic
Topology}. Clarendon Press, Oxford, 1995.

[T] C. Teleman, The quantization conjecture revisited. {\it Preprint},
math.AG/9808029.

[TZ1] Y. Tian and W. Zhang, An analytic proof of the geometric
quantization conjecture of Guillemin-Sternberg.
{\it Invent. Math.} 132 (1998), 229-259.

[TZ2] Y. Tian and W. Zhang, Holomorphic Morse inequalities in singular
reduction. {\it Math. Res. Lett.} 5 (1998), 345-352.

[TZ3] Y. Tian and W. Zhang, Symplectic reduction and a weighted multiplicity formula
for twisted Spin$^c$-Dirac operators. {\it Asian J. Math.} To appear.

[W] S. Wu, A note on higher cohomology groups of K\"ahler quotients. {\it Preprint},
math.SG/9809192.

[WZ] S. Wu and W. Zhang, Equivariant holomorphic Morse inequalities III:
non-isolated fixed points. {\it Geom. Funct. Anal.} 8 (1998), 149-178.}

$\ $

Nankai Institute of Mathematics, Nankai University, 
Tianjin 300071,
P. R. China

{\it e-mail address:} weiping@sun.nankai.edu.cn
\end{document}